\documentclass[11pt]{amsart}

\usepackage{amsmath}
\usepackage{amssymb}
\usepackage{xspace}
\DeclareMathOperator{\image}{``}
\DeclareMathOperator{\crit}{crit}
\DeclareMathOperator{\otp}{otp}
\newtheorem{theorem}{Theorem}[section]
\newtheorem{claim}[theorem]{Claim}
\newtheorem{conjecture}[theorem]{Conjecture}

\newtheorem{lemma}[theorem]{Lemma}

\newtheorem{corollary}[theorem]{Corollary}

\newcommand{\Levy}{L\'{e}vy\xspace}
\theoremstyle{definition}
\newtheorem{definition}[theorem]{Definition}

\newtheorem{question}[theorem]{Question}

\theoremstyle{remark}
\newtheorem{remark}[theorem]{Remark}

\newcount\skewfactor
\def\mathunderaccent#1#2 {\let\theaccent#1\skewfactor#2
\mathpalette\putaccentunder}
\def\putaccentunder#1#2{\oalign{$#1#2$\crcr\hidewidth
\vbox to.2ex{\hbox{$#1\skew\skewfactor\theaccent{}$}\vss}\hidewidth}}
\def\name{\mathunderaccent\tilde-3 }

\def\smallbox#1{\leavevmode\thinspace\hbox{\vrule\vtop{\vbox
   {\hrule\kern1pt\hbox{\vphantom{\tt/}\thinspace{\tt#1}\thinspace}}
   \kern1pt\hrule}\vrule}\thinspace}

\newcommand{\cf}{{\rm cf}}

\def\qedref#1{$\qed_{\reforiginal{#1}}$}

\title{The first omitting cardinal for Magidority}
\author{Shimon Garti}
\address{Institute of Mathematics,
 The Hebrew University of Jerusalem,
 Jerusalem 91904, Israel}
\email{shimon.garty@mail.huji.ac.il}

\author{Yair Hayut}
\address{School of Mathematical Sciences.
Tel Aviv University.
Tel Aviv 69978,
Israel}
\email{yair.hayut@mail.tau.ac.il}


\begin{document}
\let\labeloriginal\label
\let\reforiginal\ref

\begin{abstract}
An infinite cardinal $\lambda$ is Magidor iff $\lambda \rightarrow [\lambda]^{\aleph_0-\mathrm{bd}}_\lambda$. It is known that if $\lambda$ is Magidor then $\lambda \rightarrow [\lambda]^{\aleph_0\mathrm{-bd}}_\alpha$ for some $\alpha < \lambda$, and the first such $\alpha$ is denoted by $\alpha_M(\lambda)$. In this paper we try to understand some of the properties of $\alpha_M(\lambda)$. 
We prove that $\alpha_M(\lambda)$ can be successor of a supercompact cardinal, when $\lambda$ is a Magidor cardinal. From this result we obtain the consistency of $\alpha_M(\lambda)$ being a successor of a singular cardinal with uncountable cofinality.
\end{abstract}

\maketitle

\newpage

\section{Introduction}

The combinatorial definition of a \emph{Magidor cardinal} $\lambda$ is given by $\lambda\rightarrow[\lambda]^{\aleph_0\text{-bd}}_\lambda$. 
Recall that $[\lambda]^{\aleph_0\text{-bd}}$ is the family of all countable bounded subsets of $\lambda$. It contains subsets of any countable order type. The relation $\lambda\rightarrow[\lambda]^{\aleph_0\text{-bd}}_\lambda$ 
means that for every $c:[\lambda]^{\aleph_0\text{-bd}} \rightarrow\lambda$ there exists $A\in[\lambda]^\lambda$ for which $c\image[A]^{\aleph_0\text{-bd}}\neq\lambda$.
The reason of concentrating on bounded subsets of $\lambda$ is that $\lambda\nrightarrow[\lambda]^{\aleph_0}_\lambda$ holds for every infinite cardinal $\lambda$ (by a theorem of Erd\H{o}s and Hajnal). It follows that if $\lambda$ is Magidor then $\lambda$ is a limit cardinal of countable cofinality. 
Magidor cardinals were defined in \cite{MR3750266} through this combinatorial property.
A model-theoretic characterization of these cardinals via elementary embeddings appears in \cite{MR3666820}. 

It is easy to see that if $\lambda\rightarrow[\lambda]^{\aleph_0\text{-bd}}_\lambda$ then $\lambda\rightarrow[\lambda]^{\aleph_0\text{-bd}}_\alpha$ for some $\alpha<\lambda$ (\cite{MR3750266}, Lemma 1.2). Given a Magidor cardinal $\lambda$, the first such $\alpha$ is denoted by $\alpha_M(\lambda)$ or just $\alpha_M$ if $\lambda$ is clear from the context.
Notice that if $\beta\geq\alpha_M$ then $\lambda\rightarrow[\lambda]^{\aleph_0\text{-bd}}_\beta$ as well.
We shall say that $\alpha_M$ is \emph{the first omitting cardinal} for $\lambda$.

A parallel notion arises for J\'onsson cardinals. Recall that $\lambda$ is J\'onsson iff $\lambda\rightarrow[\lambda]^{<\omega}_\lambda$, in which case there exists a first ordinal $\alpha$ for which $\lambda\rightarrow[\lambda]^{<\omega}_\alpha$. The first such ordinal is denoted by $\alpha_J$ (or $\alpha_J(\lambda)$), and it is a regular cardinal strictly below $\lambda$.

Several properties of $\alpha_M$ are phrased in \cite{MR3750266}, among them the fact that it is always a regular cardinal. Some open problems concerning $\alpha_M$ appear in \cite{MR3750266}. A pair of related problems is labeled there as Question 3.12 and Question 3.13. The first one is whether $\alpha_M$ can be a successor of a singular cardinal. In the second one we ask about the possibility that $\alpha_M$ is a large cardinal (e.g., measurable) or the successor of a large cardinal. 
We shall prove that $\alpha_M$ is not a large cardinal but it can be a successor of a large cardinal. Thence, one can singularize this large cardinal and force $\alpha_M$ to be a successor of a singular cardinal.

Although we can force $\alpha_M$ to be a successor of a singular cardinal, the cofinality of this cardinal is uncountable in our model.
The last stage of singularizing our large cardinal is done with Magidor forcing and not with Prikry forcing. 
It is an amusing historical coincidence (Magidor cardinals were defined with no connection to Magidor forcing), but it seems that Prikry forcing fails to place $\alpha_M$ at the first point above a singular cardinal.
The reason will be explicated in the sequel, and it points to a substantial property of $\alpha_M$.

Our notation is mostly standard. We suggest \cite{MR1994835} for a comprehensive treatment to large cardinals. We shall use the Jerusalem notation in forcing, i.e. $p\leq q$ means that $q$ is stronger than $p$. 
If $\mathbb{P}$ is a forcing notion and $p,q\in\mathbb{P}$ then $p\parallel q$ means that $p$ and $q$ are compatible.
The pure order in Prikry type forcing notions will be denoted by $\leq^*$.

We call $\kappa$ supercompact iff for every ordinal $\gamma$ there exists an elementary embedding $\jmath:{\rm V}\rightarrow M$ for which $\kappa=\crit(\jmath)$ and ${}^\gamma M\subseteq M$. A forcing notion $\mathbb{P}$ is $\kappa$-directed-closed iff whenever $A\subseteq\mathbb{P}$ is a directed set of conditions and $|A|<\kappa$, there exists some $q\in\mathbb{P}$ so that $p\in A\Rightarrow p\leq q$.
We shall use Prikry and Magidor forcing, and we follow the conventions of \cite{MR2768695}. This has a particular importance with respect to Magidor forcing, as this type of forcing notions can be written in several ways.

Laver proved in \cite{MR0472529} that if $\kappa$ is supercompact then one can define a forcing notion $\mathbb{P}$ which makes $\kappa$ indestructible upon any further extension by $\kappa$-directed-closed forcing notions. The forcing $\mathbb{P}$ is a set forcing, based on the so-called Laver's diamond.

Arrows notation in this paper is coherent with the common literature. 
The notation $[\lambda]^{\aleph_0\text{-bd}}$ refers to all bounded subsets of $\lambda$ whose cardinality is $\aleph_0$, regardless of order type. We shall use the notation $[\lambda]^{\omega\text{-bd}}$ if we restrict ourselves to bounded sets of order type $\omega$.
We mention here only the less frequent relation $\lambda\rightarrow[\lambda]^{\aleph_0\text{-bd}}_{\nu,<\nu}$ which means that for every $f:[\lambda]^{\aleph_0\text{-bd}}\rightarrow\nu$ there exists $A\in[\lambda]^\lambda$ for which $|f\image[A]^{\aleph_0\text{-bd}}|<\nu$. In general, the arrows notation is designed in order to keep monotonicity, but this need not hold for $\lambda\rightarrow[\lambda]^{\aleph_0\text{-bd}}_{\nu,<\nu}$ with respect to the subscript.

We shall mention cardinals in the family of rank-into-rank, so let us recall the definitions of I0 and I1. A cardinal $\lambda$ is I1 iff there exists a non-trivial elementary embedding $\jmath:V_{\lambda+1}\rightarrow V_{\lambda+1}$. If $\lambda$ is I1 then there are many elementary embeddings of the form $k\colon V_{\lambda+1}\to V_{\lambda+1}$ with different critical points. Lest it is important to specify the critical point, we write $\mathrm{I1}(\kappa,\lambda)$ meaning that there exists an elementary embedding $\jmath\colon V_{\lambda+1}\to V_{\lambda+1}$ such that $\crit(\jmath) = \kappa$.

A cardinal $\lambda$ is I0 iff there is $\jmath:L(V_{\lambda+1}) \rightarrow L(V_{\lambda+1})$ so that $\crit(\jmath)<\lambda$. Every I1 cardinal $\lambda$ is Magidor, by a simple observation of Menachem Magidor (\cite{MR1994835}, Question 24.1).

Finally, we shall have to know that under some circumstances the assertion I1$(\kappa,\lambda)$ is preserved by forcing. We shall use, for this end, Silver's criterion. We phrase the pertinent theorem in a bit more generality than we actually need, see \cite{MR2768691}, Proposition 9.1:

\begin{theorem}
\label{ssilver} Silver's criterion. \newline 
Let $\jmath:M\rightarrow N$ be an elementary embedding, where $M,N$ are transitive models of ZFC. Let $\mathbb{P}\in M$ be a forcing notion, and $G\subseteq\mathbb{P}$ a generic subset over $M$. Assume $H$ is a generic subset of $\jmath(\mathbb{P})$ over $N$. Then the following two statements are equivalent:
\begin{enumerate}
\item [$(a)$] $\jmath(p)\in H$ for every $p\in G$.
\item [$(b)$] There exists an elementary embedding $\jmath^+:M[G]\rightarrow N[H]$ so that $\jmath^+(G)=H$ and $\jmath^+\upharpoonright M=\jmath$.
\end{enumerate}
\end{theorem}

\hfill \qedref{ssilver}

\newpage

\section{Large cardinals and their successors}

In this section we focus on Question 3.13 from \cite{MR3750266}.
Our first theorem says that under some mild assumptions, one can show that $\alpha_M$ must be a successor cardinal.
The metamathematical idea is that $\alpha_M$ (and similarly, $\alpha_J$ and $\alpha_R$ which are the parallel notions for J\'onsson and Rowbottom cardinals respectively) is not a large cardinal in the philosophical sense.

Ahead of the proof we need a lemma, which is the parallel to a lemma of Tryba, \cite{MR826499}. The lemma of Tryba refers to J\'onsson cardinals, and here we translate it to Magidor cardinals. 

\begin{lemma}
\label{ttryba} Assume that:
\begin{enumerate}
\item [$(a)$] $\lambda\rightarrow[\lambda]^{\aleph_0\text{-bd}}_{\nu,<\nu}$ and $\lambda\rightarrow[\lambda]^{\aleph_0\text{-bd}}_{\nu^{\cf(\nu)},<\nu^{\cf(\nu)}}$.
\item [$(b)$] $\nu<\lambda$ is a limit cardinal.
\item [$(c)$] There are no Magidor cardinals in the interval $[\nu,\nu^{\cf(\nu)}]$.
\end{enumerate}
Then there exists some $\rho<\nu$ for which $\lambda\rightarrow[\lambda]^{\aleph_0\text{-bd}}_{\nu,<\rho}$, and hence $\lambda\rightarrow[\lambda]^{\aleph_0\text{-bd}}_{\rho,<\rho}$.
\end{lemma}

\par\noindent\emph{Proof}. \newline 
First we show that if $\lambda\rightarrow[\lambda]^{\aleph_0\text{-bd}}_\gamma, \gamma\leq\delta\leq\varepsilon$ and there is no Magidor cardinal in the interval $[\delta,\varepsilon]$ then $\lambda\rightarrow[\lambda]^{\aleph_0\text{-bd}}_{\varepsilon,<\delta}$. Toward showing this, assume that $\lambda\nrightarrow[\lambda]^{\aleph_0\text{-bd}}_{\varepsilon,<\delta}$ and choose a function $f:[\lambda]^{\aleph_0\text{-bd}}\rightarrow\varepsilon$ so that $|f\image[A]^{\aleph_0\text{-bd}}|\geq\delta$ whenever $A\in[\lambda]^\lambda$. We may assume that $\varepsilon$ is the first counterexample.

By our assumption, $\varepsilon$ is not a Magidor cardinal. Hence $\lambda\rightarrow[\lambda]^{\aleph_0\text{-bd}}_{\varepsilon,<\varepsilon}$ (see Proposition 3.18 in \cite{MR3750266}). Let us choose $B\in[\lambda]^\lambda$ so that $\eta=|f\image[B]^{\aleph_0\text{-bd}}|<\varepsilon$. By the firsthood of $\varepsilon$ and the fact that $\delta\leq\eta<\varepsilon$ we conclude that $\lambda\rightarrow[\lambda]^{\aleph_0\text{-bd}}_{\eta,<\delta}$. In particular, one can choose a subset $C\in[B]^\lambda$ for which $|f\image[C]^{\aleph_0\text{-bd}}|<\delta$, a contradiction to the choice of $f$.

We proceed to the assertion of the lemma.
Let $A$ be the set $\{\sigma<\nu: \lambda\nrightarrow[\lambda]^{\aleph_0\text{-bd}}_{\nu,<\sigma}\}$. Denote $\sup(A)$ by $\eta$. If $\eta<\nu$ then $\eta^+<\nu$ as well (recall that $\nu$ is a limit cardinal) so $\eta^+$ can serve as the alleged $\rho$ in the lemma. 

Assume towards contradiction that $\sup(A)=\nu$, and choose a sequence of members of $A$ of the form $\langle\sigma_\alpha:\alpha<\cf(\nu)\rangle$, cofinal in $\nu$. For every $\alpha<\cf(\nu)$ choose a function $f_\alpha:[\lambda]^{\aleph_0\text{-bd}}\rightarrow\nu$ such that $|f_\alpha\image[x]^{\aleph_0\text{-bd}}|\geq\sigma_\alpha$ whenever $x\in[\lambda]^\lambda$.

Let $B=\prod\limits_{\alpha<\cf(\nu)}\sigma_\alpha$. We define a function $g:[\lambda]^{\aleph_0\text{-bd}}\rightarrow B$ as follows:
$$
g(s)=(f_\alpha(s):\alpha<\cf(\nu)).
$$
Notice that $|B|=\nu^{\cf(\nu)}$. By assumption $(a)$, 
$\lambda\rightarrow[\lambda]^{\aleph_0\text{-bd}}_{\nu^{\cf(\nu)},<\nu^{\cf(\nu)}}$.
Hence there exists a set $x\in[\lambda]^\lambda$ for which $|g\image[x]^{\aleph_0\text{-bd}}|<\nu$. This follows from the beginning of the proof, by letting $\gamma=\delta=\nu, \varepsilon=\nu^{\cf(\nu)}$, upon noticing that there are no Magidor cardinals in the interval $[\delta,\varepsilon]$ by assumption $(e)$. Consequently, $\lambda\rightarrow[\lambda]^{\aleph_0\text{-bd}}_{\nu^{\cf(\nu)},<\nu}$, which amounts to the existence of $x\in[\lambda]^\lambda$ so that $|g\image[x]^{\aleph_0\text{-bd}}|<\nu$.

On the other hand, every value of $f_\alpha$, for each $\alpha<\cf(\nu)$, gives rise to a distinct element of $B$. Hence $|g\image[x]^{\aleph_0\text{-bd}}|\geq \bigcup\limits_{\alpha<\cf(\nu)}|f_\alpha\image[x]^{\aleph_0\text{-bd}}|=\nu$, and this contradiction gives the desired conclusion.

\hfill \qedref{ttryba}

Based on this lemma, we can prove the following:

\begin{theorem}
\label{t3} Let $\lambda$ be a Magidor cardinal.
\begin{enumerate}
\item [$(a)$] If there is no Magidor cardinal in the interval $[\alpha_M,2^{\alpha_M}]$ then $\alpha_M$ is a successor cardinal.
\item [$(b)$] If every limit cardinal is a strong limit cardinal then $\alpha_M(\lambda)$ is a successor cardinal for every Magidor cardinal $\lambda$.
\end{enumerate}

\end{theorem}

\par\noindent\emph{Proof}. \newline 
As mentioned in the introduction, if $\lambda$ is Magidor then $\lambda$ is a limit cardinal.

Part $(b)$ follows from part $(a)$ by noticing that if every limit cardinal is strong limit then there are no limit cardinals in the interval $[\alpha_M,2^{\alpha_M}]$ and hence no Magidor cardinals in this interval. We prove, therefore, part $(a)$.

Assume towards contradiction that $\alpha_M$ is a limit cardinal. All the requirements of Lemma \ref{ttryba} hold, bearing in mind that $\alpha_M$ here stands for $\nu$ there. Requirement $(a)$ there is a simple property of $\alpha_M$ as proved in \cite[Theorem 1.8]{MR3750266} and $(b)$ is our assumption towards contradiction. Requirements $(c)$ is the assumption of the theorem.

It follows from the conclusion of Lemma \ref{ttryba} that $\lambda\rightarrow[\lambda]^{\aleph_0\text{-bd}}_{\rho,<\rho}$ for some $\rho<\alpha_M$, but this is an absurd in the light of the definition of $\alpha_M$, so we are done.

\hfill \qedref{t3}

Our next goal is to show that $\alpha_M$ can be basically a successor of every regular cardinal. This is possible even if one wishes to force $\alpha_M$ to be successor of large cardinals.
The following preservation theorem is in the spirit of the celebrated \Levy-Solovay preservation theorem from \cite{MR0224458} for measurable cardinals. 

\begin{claim}
\label{llevysolovay} Let $\lambda$ be Magidor and $\alpha<\lambda$. Let $\mathbb{P}$ be an $\alpha$-cc $\aleph_1$-complete forcing notion. \newline 
Then $\lambda$ remains Magidor in the generic extension by $\mathbb{P}$.
\end{claim}

\par\noindent\emph{Proof}. \newline 
Without loss of generality, $\alpha$ is regular and hence not Magidor (we can always work with $\alpha^+$ in lieu of $\alpha$). We may also assume that $\alpha_M \leq \alpha$, by taking larger $\alpha$ if needed. Let $\name{f}$ be a name of a function from $[\lambda]^{\aleph_0\text{-bd}}$ into $\alpha$, and let $p$ be a condition which forces this fact.

We define $g\in{\rm V}$, which is also a function from $[\lambda]^{\aleph_0\text{-bd}}$ into $\alpha$. Given any $t\in[\lambda]^{\aleph_0\text{-bd}}$ let $g(t)=\sup \{\eta<\alpha: \exists q, p\leq q, q\Vdash\name{f}(t)=\eta\}$. By the chain condition and the regularity of $\alpha$, $g(t)<\alpha$ for every $t\in[\lambda]^{\aleph_0\text{-bd}}$. By the assumption of $\aleph_1$-completeness, $[\lambda]^{\aleph_0\text{-bd}}$ is the same mathematical object both in ${\rm V}$ and ${\rm V}[G]$.

Choose $A\in[\lambda]^\lambda$ for which $|g\image[A]^{\aleph_0\text{-bd}}|<\alpha$. This can be done since $\alpha$ is not a Magidor cardinal. But now $\sup (\name{f}\image[A]^{\aleph_0\text{-bd}})\leq \sup (g\image[A]^{\aleph_0\text{-bd}}) <\alpha$. In particular, $\name{f}$ omits colors on a full size subset of $\lambda$, so $\lambda$ is Magidor.

\hfill \qedref{llevysolovay}

We shall see, below, that if $\lambda$ is I1 and $\kappa<\lambda$ where $\kappa$ is measurable then $\lambda$ remains Magidor after adding a Prikry sequence to $\kappa$.
Despite the possible preservation of Magidority by Prikry forcing, it turns out that a ``small" Prikry forcing may change the value of $\alpha_M$, in an interesting way. If $\lambda$ is Magidor then one can force $\alpha_M(\lambda)=\aleph_2$ while preserving the Magidority of $\lambda$ (Proposition 1.10 of \cite{MR3750266}). 
This is done, essentially, by collapsing the predecessor (or predecessors) of $\alpha_M$, and it gives a simple way to decrease $\alpha_M$.
Using Prikry forcing, one can \emph{increase} $\alpha_M$.

\begin{theorem}
\label{increasealpham} Let $\lambda$ be Magidor, and let $\kappa<\lambda$ be a measurable cardinal so that $2^\kappa<\lambda$. Let $\mathbb{P}$ be Prikry forcing through some normal ultrafilter $\mathcal{U}$ over $\kappa$. Let $G\subseteq\mathbb{P}$ be generic. \newline 
If $\lambda$ is still Magidor in ${\rm V}[G]$ then $\alpha_M>\kappa$.
Moreover, $\alpha_M>(\kappa^\omega)^{V[G]}$, so $\alpha_M>\kappa^+$ in ${\rm V}[G]$.
\end{theorem}

\par\noindent\emph{Proof}. \newline 
First we prove that if $\mu=\cf(\mu)>2^\kappa$ and $T\in[\mu]^\mu\cap{\rm V}[G]$ then there exists $S\subseteq T$ so that $S\in[\mu]^\mu\cap{\rm V}$.
For this end, assume that $\name{y}$ is a name of a subset of $\mu$ of size $\mu$, and recall that a generic subset $G$ has been chosen.
The interpretation of $\name{y}$ according to $G$ can be written as ${\name{y}}_G = \bigcup_{p\in G}y_p$ where $y_p = \{\alpha\in\mu: p\Vdash\check{\alpha}\in\name{y}\}$.
Observe that each $y_p$ belongs to the ground model, as the forcing relation is definable in ${\rm V}$.
Since $\cf({\name{y}}_G)>2^\kappa=|G|$ we see that there exists a single condition $p\in G$ and a set $y_p\in[\mu]^\mu\cap{\rm V}$ such that $p\Vdash y_p\subseteq {\name{y}}_G$, as desired.

Our objective is to define, in ${\rm V}[G]$, a function $f$ from $[\lambda]^{\aleph_0\text{-bd}}$ into $\kappa$ which omits no color on full size subsets of $\lambda$. The main point is to take care of new sets of size $\lambda$, added by the forcing poset. As a preliminary, for every $\alpha<\lambda$ such that $\cf^{\rm V}(\alpha)=\kappa$ we choose in ${\rm V}$ a cofinal sequence $\langle\beta^\alpha_i:i<\kappa\rangle$. We also fix a function $g:[\kappa]^\omega\rightarrow\kappa$, now in ${\rm V}[G]$, such that $g\image[H]^\omega=\kappa$ whenever $H\in[\kappa]^\kappa$. We may assume that $g$ is defined only over unbounded subsets of $\kappa$ (recall that $\cf^{\rm V[G]}(\kappa)=\omega$). The existence of $g$ can be proved as the existence proof of the usual $\omega$-J\'onsson functions.

Assume now that $t=\{t_n:n\in\omega\}$ belongs to $[\lambda]^{\aleph_0\text{-bd}}$. Let $\gamma_t$ be $\sup(\{t_n + 1 \mid n < \omega\})$. If $\cf^{\rm V}(\gamma_t)\neq\kappa$ then we define $f(t)=0$. Assume that $\cf^{\rm V}(\gamma_t)=\kappa$. For every $n\in\omega$ let $\rho_n$ be the first ordinal $i<\kappa$ so that $t_n\leq\beta^{\gamma_t}_i$. We define $f(t)=g(\{\rho_n:n\in\omega\})$.

Assume that $T\in[\lambda]^\lambda\cap{\rm V}[G]$. 
Choose any regular cardinal $\theta<\lambda$ such that $\theta>2^\kappa$.
By the assertion from the beginning of the proof we choose $S\in[\lambda]^\theta\cap{\rm V}$ such that $S\subseteq T$. 
Let $\gamma$ be the supremum the first $\kappa$ elements of $S$. We shall prove that $f\image[S\cap\gamma]^\omega=\kappa$, thus accomplishing the proof (notice that all the members of $[S\cap\gamma]^\omega$ are bounded in $\lambda$).

Suppose $\eta<\kappa$ is any color. Since $\langle\beta^\gamma_i:i<\kappa\rangle$ is cofinal in $\gamma$ and since $\kappa$ is regular in $V$, the set $W=\{\rho<\kappa:\exists\delta\in S, \rho=\min\{j<\kappa:\delta\leq\beta^\gamma_j\}\}$ is of size $\kappa$. By the nature of $g$, we can choose $\{\rho_n:n\in\omega\}\subseteq W$ for which $g(\{\rho_n:n\in\omega\})=\eta$. Notice that $\sup \{\rho_n \mid n < \omega\} = \kappa$. 
For each $n\in\omega$ choose $t_n\in S$ such that $\rho_n=\min\{j<\kappa: t_n\leq\beta^\gamma_j\}$, and let $t=\{t_n:n\in\omega\}$.
Now $f(t)=g(\{\rho_n:n\in\omega\})=\eta$, so we are done.

We show now how to modify the proof in order to get $\Vdash_{\mathbb{P}} \alpha_M>(\kappa^\omega)^{V[G]}$. We fix the sequences $\langle\beta^\alpha_i:i<\kappa\rangle$ as before, and the $\omega$-J\'onsson function $g$ in $V[G]$. Our goal is to define $f:[\lambda]^{\aleph_0\text{-bd}} \rightarrow \kappa^\omega$ which omits no sequence in $\kappa^\omega$ over any full size subset of $\lambda$.

Suppose that $t\in[\lambda]^{\aleph_0\text{-bd}}$. If $\otp(t)\neq\omega\cdot\omega$ then let $f(t)$ be $\vec{0}$, the fixed sequence of zeros. Likewise, if $\gamma_t=\sup(t)$ is not an ordinal of cofinality $\kappa$ in the ground model then we let $f(t) = \vec{0}$. If $\otp(t)=\omega\cdot\omega$ and $\cf^{\rm V}(\gamma_t)=\kappa$ enumerate the ordinals of $t$ by $\{\langle t_{\omega\cdot m+n}:n\in\omega\rangle :m\in\omega\}$ in increasing order. For each $m\in\omega$ denote the slice $\langle t_{\omega\cdot m+n}:n\in\omega\rangle$ by $t^m$.
For every $m,n\in\omega$ let $\rho^m_n$ be the first ordinal $i<\kappa$ such that $t^m_n\leq\beta_i^{\gamma_t}$. finally, define $f(t)=\langle g(\{\rho^m_n:n\in\omega\}):m\in\omega\rangle$.

Assume now that $T\in[\lambda]^\lambda\cap{\rm V}[G]$, and chose a sufficiently large regular $\theta < \lambda$ such that for some $S \in [\lambda]^\theta \cap V$ and some $p \in G$, we have $p \Vdash \check{S} \subseteq \name{T}$. 
We step up a bit further, and concentrate on the supremum of the first $\kappa \cdot \omega$ elements of $S$, say $\gamma$. We claim that $f\image[S\cap\gamma]^{\aleph_0}=\kappa^\omega$.

Assume that $\langle\eta_m:m\in\omega\rangle\in\kappa^\omega$. Let $\gamma_{-1}$ be $0$ and for every $m\in\omega$ let $\gamma_m$ be the supremum of the first $\kappa$ elements of $S$ above $\gamma_{m-1}$. For every $m\in\omega$ let \[W_m=\{\rho<\kappa: \exists\delta\in S,\gamma_{m-1}\leq\delta<\gamma_{m},\ \rho=\min\{j<\kappa: \delta\leq\beta^{\gamma_m}_j\}\}.\] By the choice of $g$ we choose, for each $m\in\omega$, a sequence $(\rho^m_n:n\in\omega)\subseteq W_m$ for which $g(\{\rho^m_n:n\in\omega\})=\eta_m$.

Now we can finish the proof as follows. For every $m\in\omega$ let $t^m=\{t^m_n:n\in\omega\}$ where $t^m_n\in S$ satisfies $t^m_n\leq\beta^{\gamma_m}_j$ with respect to $j=\rho^m_n$ (and assuming that $\rho^m_n$ is the first $j$ with this property). Define $t=\bigcup\limits_{m\in\omega}t^m$ and observe that $f(t)=\langle g(t^m):m\in\omega\rangle = \langle \eta_m:m\in\omega\rangle$, thus accomplishing the proof.

\hfill \qedref{increasealpham}

The fact that Prikry forcing through $\kappa$ results in $\alpha_M\geq\kappa^{++}$ is the strongest reason which stands behind Conjecture \ref{conj0} below. Namely, we suspect that $\alpha_M$ cannot be a successor of a singular cardinal with countable cofinality. The most natural way to force it is Prikry forcing, and this provably fails.
An analysis of the proof shows that specific properties of Prikry forcing are not essential for the validity of the basic argument.
The main point in the above proof can be abstracted as follows.

\begin{corollary}
\label{corcov} Let $V,W$ be models of ZFC. \newline 
Assume that $V\subseteq W$ and $\lambda$ is Magidor in both of them.
Assume further that $\mu<\lambda, \mu>\cf(\mu)=\omega$ in $W$ and $\mu$ is regular in $V$. \newline 
If every $S\in[\lambda]^\lambda\cap W$ contains a set $T\in V$ of order type $\mu \cdot \omega$, then $\alpha_M>\mu^+$ in $W$.
\end{corollary}

\hfill \qedref{corcov}

If $\lambda$ is Magidor and $\kappa<\lambda$ then $\kappa^+<\lambda$, so the chain condition of Prikry forcing through $\kappa$ is promising. However, Prikry forcing is not $\aleph_1$-complete, so it may ruin the Magidority of $\lambda$ or change the value of $\alpha_M$. In order to employ Prikry forcing we need stronger assumptions.
By I1$(\kappa,\lambda)$ we mean that $\lambda$ is I1 as witnessed by $\jmath: {\rm V}_{\lambda+1}\rightarrow {\rm V}_{\lambda+1}$ and $\kappa=\crit(\jmath)$. If $\mu$ is a measurable cardinal below $\kappa$ then Prikry forcing for $\mu$ preserves I1$(\kappa,\lambda)$. More generally, any small forcing notion keeps I1. 
Probably, the following Lemma is known, but we elaborate:

\begin{lemma}
\label{ppp} \Levy-Solovay for I1. \newline 
Assume $I1(\kappa,\lambda),$ and $\mathbb{P}\in V_\kappa$ is a forcing notion. Then $I1(\kappa,\lambda)$ holds in ${\rm V}^{\mathbb{P}}$.
\end{lemma}

\par\noindent\emph{Proof}. \newline 
Assume $\jmath:V_{\lambda+1}\rightarrow V_{\lambda+1}$ witnesses $I1(\kappa,\lambda)$. Since $\mathbb{P}\in V_\kappa$ we know that $\jmath(p)=p$ for every $p\in\mathbb{P}$, and likewise $\jmath(\mathbb{P})=\mathbb{P}$. Fix a generic subset $G\subseteq\mathbb{P}$. We will use Silver's criterion with $M=N=V_{\lambda+1}$ and $H=G$. Since our formulation of Theorem \ref{ssilver} does not immediately apply to this case, we will continue and give a detailed proof.

We claim that there is an elementary embedding $\jmath^+$ from $V_{\lambda+1}^{V[G]}$ into $V_{\lambda+1}^{V[G]}$ which extends $\jmath$. This implies, in particular, that $\kappa=\crit(\jmath^+)$ and hence $V[G]\models I1(\kappa,\lambda)$. For proving this claim notice that if $\name{x}$ is a name of an element in $V_{\lambda+1}\cap V[G]$ then for some name $\name{y}\in V_{\lambda+1}$ we have $\Vdash_{\mathbb{P}} \name{x}=\name{y}$, so we can focus only on names which belong to $V_{\lambda+1}$.

Given a $\mathbb{P}$-name which belongs to $V_{\lambda+1}$, let $\jmath^+({\name{y}}_G)=(\jmath(\name{y}))_G$. If $\check{y}$ is a canonical name of an element $y$ in $V_{\lambda+1}^{\rm V}$ then $\jmath^+(y)= \jmath^+(\check{y}_G)=(\jmath(\check{y}))_G=\jmath(y)$, the last equality follows from the elementarity of $\jmath$. We conclude that $\jmath^+$ extends $\jmath$.

Similarly, we argue that $\jmath^+$ is elementary. For this, let $\varphi$ be any first order formula and $\name{y}$ a name in $\mathbb{P}$. We see that:
$$
V_{\lambda+1}^{V[G]}\models \varphi[{\name{y}}_G] \Leftrightarrow \exists p\in G, p \Vdash \varphi[\name{y}] \Leftrightarrow 
$$

$$
\exists p\in G, p \Vdash \varphi[\jmath(\name{y})] \Leftrightarrow V_{\lambda+1}^{V[G]}\models \varphi[\jmath(\name{y})_G].
$$
By this, $\jmath^+$ is elementary in the generic extension, so we are done.

\hfill \qedref{ppp}

\begin{corollary}
\label{cori1} Assume $I1(\mu_\alpha,\lambda_\alpha)$ where $\langle\mu_\alpha: \alpha\in{\rm Ord}\rangle$ is a proper class. \newline 
Then one can force the existence of a Magidor cardinal with $\alpha_M$ arbitrarily large. Likewise, it is consistent that $\lambda$ is Magidor and the distance between $\alpha_J$ and $\alpha_M$ is arbitrarily large.
\end{corollary}

\par\noindent\emph{Proof}. \newline 
For the first assertion choose any measurable cardinal $\kappa$ and add a Prikry sequence to it. Then use Theorem \ref{increasealpham} in order to conclude that all instances of I1 with critical point above $\kappa$ are still I1 and hence Magidor, with $\alpha_M$ above $\kappa$.

For the second assertion notice that each I1 cardinal $\lambda$ is an $\omega$-limit of measurable cardinals, hence Rowbottom. It follows that for colorings of finite sets of $\lambda$ we can find a full size subset which assumes only countably many colors, i.e. $\alpha_J=\aleph_1$. Now use the former paragraph to increase $\alpha_M$ while keeping $\lambda$ as I1, so its $\alpha_J$ is still $\aleph_1$.

\hfill \qedref{cori1}

\begin{remark}
\label{r1} The general fact proved above, that if $\mu=\cf(\mu)>2^\kappa$ and $T\in[\mu]^\mu\cap{\rm V}[G]$ then there exists $S\subseteq T$ so that $S\in[\mu]^\mu\cap{\rm V}$, shows that a small Prikry forcing cannot \emph{create} a new Magidor cardinal. Namely, if $\lambda>\cf(\lambda)=\omega$ is not Magidor, and $\kappa$ is a measurable cardinal for which $2^\kappa<\lambda$ then Prikry forcing through $\kappa$ keeps the non-Magidority of $\lambda$. 
This should be compared with a theorem of Woodin, \cite{MR2914848}, who proved that if I0$(\kappa,\lambda)$ then Prikry forcing through $\kappa$ makes $\kappa$ I1, and hence Magidor.
\end{remark}

\hfill \qedref{r1}

Merging the above method with \Levy collapses, we can show that basically $\alpha_M$ can be any prescribed successor of a regular cardinal. Moreover, it is consistent that $\alpha_M$ is a successor of a strongly inaccessible cardinal or a strongly Mahlo cardinal. For proving this, we need another lemma about the impact of \Levy collapse on $\alpha_M$.

\begin{lemma}
\label{qqq} Collapsing $\alpha_M$. \newline 
Assume $\aleph_0 < \mu=\cf(\mu)<\lambda, \lambda$ is Magidor and $\mu^+<\alpha_M(\lambda)\leq\alpha = \cf(\alpha)<\lambda$. 
Assume further that $\alpha$ is $\mu$-closed (i.e.\ for all $\beta < \alpha$, $\beta^\mu < \alpha$).
Let $\mathbb{P}=L\acute{e}vy(\mu,<\alpha)$ and let $G\subseteq\mathbb{P}$ be generic. 
Then $V[G]\models\alpha_M(\lambda)=\mu^+$.
\end{lemma}

\par\noindent\emph{Proof}. \newline
We have to prove the following two statements:
\begin{enumerate}
\item [$(\aleph)$] $\Vdash_{\mathbb{P}} \alpha_M\leq\mu^+$.
\item [$(\beth)$] $\Vdash_{\mathbb{P}} \alpha_M\geq\mu^+$.
\end{enumerate}
The first assertion follows from the chain condition. By the assumption of the lemma, $\mathbb{P}$ is $\alpha$-cc. Now let $\name{f}$ be a name and let us fix a condition $p\in\mathbb{P}$ that forces that $\name{f}$ is a function from $[\lambda]^{\aleph_0\text{-bd}}$ into $\alpha$. We define another function $g:[\lambda]^{\aleph_0\text{-bd}}\rightarrow\alpha, g\in{\rm V}$, as follows. Given $s\in[\lambda]^{\aleph_0\text{-bd}}$ let $g(s)=\sup\{\beta<\mu^+: \exists q\geq p, q\Vdash\name{f}(s)=\beta\}$.
Notice that $g(s)\in\alpha$ since $\alpha$ is regular and by the chain condition.

In the ground model we have $\lambda\rightarrow[\lambda]^{\aleph_0\text{-bd}}_\alpha$, so we choose a set $A\in[\lambda]^\lambda$ for which $|g\image[A]^{\aleph_0\text{-bd}}|<\alpha$. This can be done since $\alpha$ is regular and hence not Magidor.  
Fix an ordinal $\gamma < \alpha$ so that $g\image[A]^{\aleph_0\text{-bd}}\subseteq\gamma$, and notice that $p\Vdash_{\mathbb{P}}\sup\{g(s):s\in [A]^{\aleph_0\text{-bd}}\}\geq \sup\{\name{f}(s):s\in [A]^{\aleph_0\text{-bd}}\}$, since $\mathbb{P}$ is $\aleph_1$-complete so no new countable sets are forced into the universe. We conclude that $p\Vdash_{\mathbb{P}} \name{f}\image[A]^{\aleph_0\text{-bd}}\subseteq\gamma<\alpha$, which means that $V[G]\models\lambda\rightarrow[\lambda]^{\aleph_0\text{-bd}}_\alpha$.
By the collapse, $V[G]\models\alpha=\mu^+$, so part $(a)$ is accomplished.

The second assertion follows from the size of $\mathbb{P}$. We have to find a coloring which exemplifies $\lambda\nrightarrow[\lambda]^{\aleph_0\text{-bd}}_\mu$ in the generic extension. Choose a coloring $c:[\lambda]^{\aleph_0\text{-bd}}\rightarrow\mu$ which shows that $\lambda\nrightarrow[\lambda]^{\aleph_0\text{-bd}}_\mu$ in ${\rm V}$. We claim that $\check{c}$ gives the same relation in $V[G]$.

For this recall that $[\lambda]^{\aleph_0\text{-bd}}$ is the same object in ${\rm V}$ and in $V[G]$ by $\aleph_1$-completeness. The only possible problem would be a new $\name{A}\in[\lambda]^\lambda$ which might omit colors. In order to cope with this problem we shall prove that for some $B\in[\lambda]^\lambda\cap{\rm V}$ we have $\Vdash_{\mathbb{P}} \check{B}\subseteq\name{A}$.

Indeed, for each ordinal $\beta<\lambda$ let $\sigma_\beta$ be the statement $\check{\beta}\in\name{A}$. Let $\langle\lambda_n:n\in\omega\rangle$ be an increasing sequence of regular cardinals such that $\alpha<\lambda_0$ and $\lambda=\bigcup\limits_{n\in\omega}\lambda_n$. By induction on $n\in\omega$ we shall define a set $B_n\in[\lambda_n]^{\lambda_n}$ and a condition $q_n$ such that $q_n\Vdash \check{B}_n\subseteq\name{A}$. Moreover, the sequence $\langle q_n:n\in\omega\rangle$ will be increasing.

Suppose $q_m,B_m$ were constructed for every $m<n$. For $\lambda_n$-many ordinals $\beta$ there is a condition which forced $\sigma_\beta$. Since $\lambda_n=\cf(\lambda_n)>|\mathbb{P}|$ we can pick a single condition $q_n$ (above $q_{n-1}$ in case $n>0$) such that $B_n=\{\beta<\lambda_n:q_n\Vdash\sigma_\beta\}$ is of size $\lambda_n$.

By the completeness of $\mathbb{P}$ we choose a condition $q$ so that $q\geq q_n$ for every $n\in\omega$. Let $B=\bigcup\limits_{n\in\omega}B_n$. Notice that $q\Vdash\check{B}\subseteq\name{A}$ and $|B|=\lambda$, so $c\upharpoonright [B]^{\aleph_0\text{-bd}}=\mu$. Hence $\Vdash_{\mathbb{P}} \check{c}\image [\name{A}]^{\aleph_0\text{-bd}}=\mu$, as desired.

\hfill \qedref{qqq}

For the purpose of forcing $\alpha_M$ to be a successor of a singular cardinal we shall need to force that $\alpha_M$ is a successor of a measurable cardinal, and this will be done later.
But the above claims enable us to show that $\alpha_M$ can be a successor of small large cardinals.

\begin{claim}
\label{c1} Making $\alpha_M$ successor of small large cardinals.
\begin{enumerate}
\item [$(a)$] For every successor ordinal $\beta$, it is consistent (assuming the existence of large cardinals) that $\alpha_M(\lambda)=\aleph_{\beta+1}$ for some Magidor cardinal $\lambda$.
\item [$(b)$] It is consistent that $\alpha_M(\lambda)$ is a successor of a strongly inaccessible cardinal (and even strongly Mahlo).
\end{enumerate}
\end{claim}

\par\noindent\emph{Proof}. \newline 
We need the assumption ${\rm I1}(\kappa,\lambda)$ for some $\kappa$ above $\aleph_{\beta}$, so that there exists a measurable cardinal $\mu<\kappa$ with $\aleph_{\beta}<\mu$. Now we force with Prikry forcing through $\mu$, so if $G$ is a generic set for Prikry forcing then $\lambda$ is still Magidor in ${\rm V}[G]$ and $\alpha_M(\lambda)>\mu$.

The next stage is to collapse the predecessors of $\alpha_M$ to $\aleph_{\beta}$. In $V[G]$ choose a regular $\alpha$ such that $\alpha_M \leq \alpha < \lambda$ and $\alpha$ is $\aleph_\beta$-closed (such $\alpha$ exists, since $\lambda$ is still a limit of inaccessible cardinals). Let $H$ be a generic set in ${\rm V}[G]$ for the collapse $\text{\Levy}(\aleph_\beta, <\alpha)$. It follows from Lemma \ref{qqq} that $\alpha_M(\lambda)=\aleph_{\beta}^+$ in ${\rm V}[G][H]$, so we are done with part $(a)$.

For part $(b)$ notice that the collapse (being complete enough) adds no bounded subsets to the predecessor of $\alpha_M$. Hence if this is an inaccessible cardinal in ${\rm V}$ then it is still inaccessible in ${\rm V}[G][H]$. A similar argument shows that $\alpha_M$ can be forced to be a successor of a strongly Mahlo cardinal.

\hfill \qedref{c1}

Question 3.13 from \cite{MR3750266} asks for a stronger statement: Is it consistent that $\alpha_M(\lambda)$ is a successor of a measurable cardinal? We shall prove that a positive answer is consistent, even if one replaces measurability by supercompactness.
Basically, we would like to lift $\alpha_M$ above some supercompact cardinal and then to collapse its predecessors to this supercompact. Prikry forcing is a useful way to achieve the first mission, but it ruins supercompactness below it since it adds a weak square (and even stronger forms of squares, see \cite{MR1360144}). Fortunately, we also have a delicate way to increase $\alpha_M$, based on the quilshon principle from \cite{MR3750266}:

\begin{definition}
\label{ssss} Quilshon. \newline
Assume $\lambda>\delta=\cf(\delta)$. \newline
We say that $\pitchfork_{\lambda,\delta}$ holds iff there is a collection $\{S_\gamma:\gamma<\delta\}$ of disjoint subsets of $\lambda$ so that $S_\gamma\cap\eta$ is a stationary subset of $\eta$ for every ordinal $\eta<\lambda$ with $\cf(\eta)=\delta$ and every $\gamma<\delta$.
\end{definition}

It has been proved in \cite{MR3750266}, Theorem 2.2, that $\pitchfork_{\lambda,\delta}$ implies $\alpha_M(\lambda)>\delta$. Likewise, adding $\pitchfork_{\lambda,\delta}$ by the partial square forcing of Jensen preserves supercompactness below $\delta$ (Theorems 2.6 and 2.8 of \cite{MR3750266}). This yields the consistency of $\alpha_M$ being a successor of a supercompact cardinal:

\begin{theorem}
\label{mt} It is consistent that $\lambda$ is Magidor, $\alpha_M(\lambda)=\mu^+$ and $\mu$ is supercompact.
\end{theorem}

\par\noindent\emph{Proof}. \newline 
We begin with $I1(\kappa,\lambda)$ and we fix a supercompact cardinal $\mu$ below $\kappa$. This can be arranged if we choose $\lambda$ to be a limit of supercompact cardinals, since $\kappa=\crit(\jmath)$ can be arbitrarily large below $\lambda$ as the $n$-th iteration elementary embedding $\jmath^n$ satisfies $\crit(\jmath^n)=\jmath^n(\kappa)$ and the values of $\jmath^n(\kappa)$ are unbounded in $\lambda$.

Choose a regular cardinal $\delta\in(\mu,\kappa)$ and force $\pitchfork_{\lambda,\delta}$ while preserving the supercompactness of $\mu$ on the one hand and I1$(\kappa,\lambda)$ on the other hand. The canonical way to force $\pitchfork_{\lambda,\delta}$ gives these properties (see Theorems 2.6 and 2.8 of \cite{MR3750266}).
We force now with Laver's forcing, making $\mu$ indestructible upon $\mu$-directed-closed forcing notions. By virtue of Lemma \ref{ppp}, $I1(\kappa,\lambda)$ holds in the generic extension and $\mu<\alpha_M(\lambda)$. We may assume that $\alpha_M(\lambda) < \kappa$ since $I1(\kappa,\lambda)$ implies that there are arbitrary large cardinals $\kappa' < \lambda$ such that $I1(\kappa',\lambda)$ holds. By changing $\kappa$ to a larger $\kappa'$ if needed, we may arrange that $\alpha_M(\lambda) < \kappa$.

If $\alpha_M(\lambda)=\mu^+$ we are done. If not, let $\alpha=((\alpha_M)^\mu)^+$ and notice that $\alpha<\kappa$ since $\kappa$, being measurable, is strongly inaccessible.

Let $\mathbb{P}=L\acute{e}vy(\mu,<\alpha)$ and choose a generic subset $G\subseteq\mathbb{P}$. Since $\mathbb{P}$ is $\mu$-directed closed, forcing with $\mathbb{P}$ preserves the supercompactness of $\mu$.
By Lemma \ref{ppp} we still have I1$(\kappa,\lambda)$, so in particular $\lambda$ remains a Magidor cardinal. By Lemma \ref{qqq}, $V[G]\models\alpha_M(\lambda)=\mu^+$, so the proof is accomplished.

\hfill \qedref{mt}

We conclude this section with two open problems:

\begin{question}
\label{q1} Is it consistent, under any assumption, that $\alpha_M$ is a limit cardinal?
\end{question}

The second question is about $\alpha_J$. Our knowledge about $\alpha_J$ is relatively poor (see \cite{MR826499}). 
We know how to obtain a J\'onsson cardinal with large $\alpha_J$ but we do not know how to change $\alpha_J$ for a given cardinal.
The following is typical:

\begin{question}
\label{q0} Can we increase $\alpha_J$ to an arbitrarily large value?
\end{question}

\newpage 

\section{Successors of singular cardinals}

In this section we focus on Question 3.12 from \cite{MR3750266}, namely is it consistent that $\alpha_M(\lambda)=\mu^+$ when $\mu$ is a singular cardinal?
Our approach depends on the tentative answer to the problem.
If we speculate that the answer is no, then the most natural thing would be to express $\mu^+$ as ${\rm tcf}(\prod_{\alpha < \cf(\mu)}\mu_\alpha,J)$ where $\langle \mu_\alpha \mid \alpha < \cf (\mu)\rangle$ is an increasing sequence of regular cardinals that is cofinal in $\mu$, and $J$ is an ideal over $\cf(\mu)$.
Now one can fix $f_\alpha:[\lambda]^{\aleph_0\text{-bd}} \rightarrow \mu_\alpha$ which omits no colors, for every $\alpha<\cf(\mu)$.
The hope is to define from these functions a coloring $c:[\lambda]^{\aleph_0\text{-bd}} \rightarrow \mu^+$ which omits no color over full-sized subsets of $\lambda$.

If one wishes to try a positive answer then the natural attempt would be to force $\alpha_M=\mu^+$ where $\mu$ is measurable (or even more) and then to singularize $\mu$.
If this process keeps $\alpha_M=\mu^+$ then a positive answer to the above question has been given.

Practically, there are obstacles in both ways.
It seems that there is no simple way to combine functions into small cardinals in order to get a single function $c:[\lambda]^{\aleph_0\text{-bd}} \rightarrow \mu^+$.
Actually, the main theorem of this section shows that $\alpha_M(\lambda)$ can be $\mu^+$ where $\mu>\cf(\mu)>\omega$, thus proving that this approach fails in general, though it may be helpful in case of singular cardinals with countable cofinality.

The other approach is problematic as well. 
The simplest attempt to force $\alpha_M=\kappa^+$ where $\kappa$ is measurable and then to add a Prikry sequence to $\kappa$, fails.
By Theorem \ref{increasealpham}, $\alpha_M>\kappa^\omega$ in the generic extension.
Since $\kappa>\cf(\kappa)=\omega$ after Prikry forcing, $\alpha_M>\kappa^+$ in the generic extension.
Actually, we believe that this is a ZFC limitation:

\begin{conjecture}
\label{conj0} Let $\lambda$ be Magidor and $\alpha=\alpha_M(\lambda)$. \newline 
If $\theta<\alpha$ then $\theta^\omega<\alpha$. In particular, $\alpha_M$ cannot be $\mu^+$ when $\mu>\cf(\mu)=\omega$.
\end{conjecture}

\hfill \qedref{conj0}

There is, however, an alternative to Prikry forcing.
We shall use Magidor forcing in order to force $\alpha_M$ to be a successor of a singular cardinal.
As a warm-up we show that under some assumptions on the Magidor cardinal $\lambda$ one can force $\alpha_M(\lambda)=\kappa^{++}$ by singularizing a measurable cardinal $\kappa$.
This will be done with the usual Prikry forcing through $\kappa$, so $\cf(\kappa)^{V[G]}=\omega$ and yet Prikry forcing does not increase $\alpha_M$ too much. Recall:

\begin{definition}
\label{defstrong} Strong Magidority. \newline 
Assume that $\beta<\lambda$.
\begin{enumerate}
\item [$(\aleph)$] $\lambda$ is $\beta$-Magidor iff $\lambda\rightarrow [\lambda]^{<\beta\text{-bd}}_\lambda$.
\item [$(\beth)$] $\lambda$ is strongly Magidor iff $\lambda$ is $\beta$-Magidor for every $\beta<\lambda$.
\end{enumerate}
\end{definition}

It has been shown in \cite{MR3750266} that if $\lambda$ is I1 then $\lambda$ is strongly Magidor.
If $\lambda$ is strongly Magidor and $\kappa<\lambda$ (or even $\beta$-Magidor and $\kappa\leq\beta<\lambda$) then we can define $\alpha_M^{<\kappa}(\lambda)$ as the first $\alpha<\lambda$ such that $\lambda\rightarrow[\lambda]_\alpha^{<\kappa\text{-bd}}$.
The usual $\alpha_M(\lambda)$ is then $\alpha_M^{<\omega_1}(\lambda)$.

\begin{claim}
\label{clmprikry} Let $\lambda$ be $\kappa^+$-Magidor where $\kappa<\lambda$ is a measurable cardinal.
Assume that $\alpha_M^{<\kappa^+}(\lambda)=\kappa^{++}$, $\mathbb{P}$ is Prikry forcing through $\kappa$ and $G\subseteq\mathbb{P}$ is generic. \newline 
Then $V[G]\models\alpha_M(\lambda)=\kappa^{++}$, and in particular $\lambda$ is still Magidor in $V[G]$.
\end{claim}

\par\noindent\emph{Proof}. \newline 
Let $\name{f}:[\lambda]^{\aleph_0\text{-bd}}\rightarrow \kappa^{++}$ be a name of a coloring.
We have to find $A\in[\lambda]^\lambda$ such that the interpretation of $\name{f}$ restricted to countable subsets of $A$ omits colors from $\kappa^{++}$.
It follows from the arguments of \cite{MR3750266} that $2^\kappa<\alpha_M^{<\kappa^+}(\lambda)$, so under our assumption that $\alpha_M^{<\kappa^+}(\lambda)=\kappa^{++}$ we see that $2^\kappa=\kappa^+$. The proofs in \cite{MR3750266} deal with $\aleph_0$-Magidor cardinals, but the same proofs work for the more general case as well.

For every $x\in[\lambda]^{<\kappa^+\text{-bd}}\cap V$ we define $g(x)$ to be the supremum over all ordinals $\gamma < \kappa^+$ such that there is a $\mathbb{P}$-name $\name{\tau}$ and a condition $q\in\mathbb{P}$ such that $q$ forces that $\name{\tau}$ is a countable subset of $x$ and that $\name{f}(\name{\tau}) = \gamma$.

Notice that $g\in V$ as the forcing relation is definable in $V$.
Observe also that $g(x)\in\kappa^{++}$ for every $x\in[\lambda]^{<\kappa^+\text{-bd}}\cap V$.
This is true since the number of names $\name{\tau}$ (up to equivalence) which appear in the definition of $g(x)$ is at most $(\kappa \cdot 2^\kappa)^{\kappa \cdot \aleph_0} =\kappa^+$. This is true as we can always assume that the conditions in the name $\tau$ consist of countably many maximal antichains, each of size at most $\kappa$, and there are at most $\kappa$ many ordinals in $x$. 
By the chain condition of $\mathbb{P}$ the value of $\name{f}(\name{\tau})$ can be determined in at most $\kappa$ many different ways. 

As each value of $\name{f}(\name{\tau})$ is an ordinal in $\kappa^{++}$ (recall that the range of $\name{f}$ is $\kappa^{++}$) we see that $g(x)\in\kappa^{++}$.

By the assumption that $\alpha_M^{<\kappa^+}(\lambda)=\kappa^{++}$ in the ground model, there is a set $A\in[\lambda]^\lambda$ and an ordinal $\mu<\kappa^{++}$ such that ${\rm Rang}(g)\subseteq\mu$.
It follows that ${\rm Rang}(\name{f})\subseteq\mu$ in the generic extension, so we are done.

\hfill \qedref{clmprikry}

Our next goal is to show that the assumptions of the above claim are consistent.

\begin{claim}
\label{clmprepar} Assume that $\lambda$ is I1.
\begin{enumerate}
\item [$(\aleph)$] It is consistent that $\mu<\lambda$, $\mu$ is supercompact and $\alpha_M^{<\mu^+}(\lambda)=\mu^{++}$.
\item [$(\beth)$] It is consistent that $\mu<\lambda$, $\mu$ is supercompact and $\alpha_M^{<\mu}(\lambda)=\mu^+$.
\end{enumerate}
\end{claim}

\par\noindent\emph{Proof}. \newline 
Choose $\mu<\kappa<\lambda$ so that I1$(\kappa,\lambda)$ and $\mu$ is supercompact.
Notice that $\alpha_M^{<\mu^+}(\lambda)>\mu^+$ and $\alpha_M^{<\mu}(\lambda)>\mu$. 

Let $\alpha = ((\alpha_M^{<\mu^+}(\lambda))^\mu)^+$.
Since $\lambda$ is a strong limit cardinal we see that $\alpha<\lambda$.
Now we force with $\mathbb{Q} = \text{\Levy}(\mu^{+},<\alpha)$ and one can verify that $\alpha_M^{<\mu^+}(\lambda)=\mu^{++}$ in the generic extension by $\mathbb{Q}$, as done in the proof of Theorem \ref{mt}.

A similar argument proves part $(\beth)$, but here it is possible to force $\alpha_M^{<\mu}(\lambda)=\mu^+$.
Indeed, the pertinent chain condition will be just $\mu^+$-cc, as the elements that we color are of size strictly less than $\mu$.
Consequently, $\text{\Levy}(\mu,<\alpha)$ is sufficient.

\hfill \qedref{clmprepar}

The fact that Prikry forcing at $\mu$ might increase $\alpha_M$ only to $\mu^{++}$ is suggestive.
A more illuminating formulation of this fact is that basically (under some convenient assumptions that we made) Prikry forcing at $\mu$ makes $\alpha_M=(\mu^\omega)^+$.
Philosophically this is the correct point for $\alpha_M$ since one can code $\omega$-sequences of $\mu$ in the generic extension by old $\mu$-sets, thus covering all the colors of $\mu^\omega$ but maintaining $(\mu^\omega)^+$ to omit colors.

Mathematically it suggests that if we singularize $\mu$ in such a way that $\mu^\omega=\mu$ in the generic extension then we can force $\alpha_M = (\mu^\omega)^+=\mu^+$ while $\mu>\cf(\mu)$.
This is hopeless with Prikry forcing but it can be done by another Prikry-type forcing notion which makes $\mu>\cf(\mu)>\omega$.
We let Magidor forcing into the discussion at this point.

A natural question arose in the wake of Prikry's work.
Is it possible to change the cofinality of a measurable cardinal $\kappa$ into uncountable cofinality without collapsing cardinals?
A positive answer was given by Magidor in \cite{MR0465868}, nowadays known as Magidor forcing.
The definition of the forcing is more involved than the classical Prikry forcing, and in particular the required largeness of the cardinal which changes its cofinality is much more than just measurability.

There are other differences between Prikry and Magidor forcing, the most important for us is mirrored in the covering properties of countable sets.
If $\mathbb{P}$ is Prikry forcing for $\kappa$ and $x = \{x_n:n\in\omega\}$ is a cofinal Prikry sequence, then $x$ cannot be covered by a set of size less than $\kappa$ from the ground model.
But if $\mathbb{M}$ is Magidor forcing then the situation changes a bit.
The following claim shows that Magidor forcing has a better covering property when new countable sets are considered.

\begin{claim}
\label{clmcovering} Assume that $\kappa\leq\lambda$ and $\kappa$ is sufficiently large (e.g. $\kappa$ is supercompact). Let $\mathbb{M}$ be Magidor forcing which makes $\kappa>\cf(\kappa)>\omega$. \newline 
If $\name{\tau}$ is any $\mathbb{M}$-name of an element in $[\lambda]^{\aleph_0\text{-bd}}$ then there are $p\in\mathbb{M}, \theta<\kappa$ and $x\in V$ such that $|x|=\theta$ and $p\Vdash\name{\tau}\subseteq\check{x}$.
\end{claim}

\par\noindent\emph{Proof}. \newline 
For precise definitions and explanation of facts about the Magidor forcing we refer the reader to \cite[Section 5]{MR2768695}, in which the forcing is defined using measure sequences and in particular to \cite[Subsection 5.2]{MR2768695} for details about the version which is used here. 

A condition $p\in\mathbb{M}$ is a finite sequence $\langle d_1,\ldots,d_n,(\kappa,A)\rangle$ where each $d_i$ is either an ordinal or a pair $(\nu,A_\nu)$. The part $\langle d_1,\ldots,d_n\rangle$ is the stem of $p$, and if $p,q$ share the same stem then $p\parallel q$. 
Inasmuch as the number of possible stems is $\kappa$, $\mathbb{M}$ is $\kappa^+$-cc.

Fix a generic subset $G\subseteq\mathbb{M}$.
Let $\name{\tau}$ be a name of a countable set of ordinals. 
Each element in ${\name{\tau}}_G$ is determined by an antichain of size at most $\kappa$, hence by collecting all the possibilities we construct (in $V$) a set $B$ of ordinals, $|B|\leq\kappa$ such that $\Vdash_{\mathbb{M}}\name{\tau} \subseteq \check{B}$.
Since $B\in V, |B|=\kappa$ one can fix a bijection $h:\kappa\rightarrow B$.

Let $\name{\sigma}$ be a name for a countable set of ordinals of $\kappa$, such that $\Vdash_{\mathbb{M}}\name{\sigma} = h^{-1}(\name{\tau})$.
The interpretation of $\name{\sigma}$ in the generic extension is a bounded subset of $\kappa$, since $V[G]\models \cf(\kappa)>\omega$.
Hence for some $\beta\in\kappa$ there is a condition $p\in\mathbb{M}$ such that $p\Vdash\name{\sigma}\subseteq\check{\beta}$.
It follows that $p\Vdash\name{\tau}\subseteq h\image\beta$.
Denote $h\image\beta$ by $x$ and $|\beta|=|x|$ by $\theta$.
Observe that $p$ forces that $\name{\tau}$ is contained in the set $x$ and $x\in V$. Since $\theta<\kappa$, we are done.

\hfill \qedref{clmcovering}

One can modify the above claim to elements in $[\lambda]^{\eta\text{-bd}}$ where $\eta>\aleph_0$, provided that $\mathbb{M}$ forces $\cf(\kappa)>\eta$.
In this way one can obtain $\alpha_M^{<\eta}(\lambda)=\mu^+$ where $\mu>\cf(\mu)>\eta$ as we shall prove below.
We focus on the usual $\alpha_M$ (that is, with respect to countable sets).

Let us start with a technical lemma.
\begin{lemma}\label{lemma:notJonsson}
Let $\mu < \lambda$ be cardinals and let us assume that $\mu$ is regular, uncountable and $2^{<\mu} = \mu$. 

If $\lambda\rightarrow[\lambda]^{{<}\mu\text{-bd}}_{\theta}$ and $\theta$ is not J\'onsson then $\lambda\rightarrow[\lambda]^{{<}\mu\text{-bd}}_{\theta, <\theta}$.
\end{lemma}
\par\noindent\emph{Proof}. \newline 
Let us fix a function $h \colon \mu \to (P_\mu\mu)^{<\omega}$ such that for all $x \in (P_\mu\mu)^{<\omega}$ there are unboundedly many $\zeta < \mu$ such that $h(\zeta) = x$ and if $h(\alpha) = \langle x_0, \dots, x_{n-1}\rangle$ then $x_0, \dots, x_{n-1} \subseteq \alpha$.

Let us assume that $f\colon [\lambda]^{{<}\mu\text{-bd}}\to \theta$ is a function such that for all $A \in [\lambda]^\lambda$, $|f\image [A]^{{<}\mu\text{-bd}}| = \theta$. Let $g\colon \theta^{{<}\omega}\to \theta$ witness the negative partition relation $\theta\nrightarrow [\theta]^{{<}\omega}_\theta$. 

Let us define a function $F\colon [\lambda]^{{<}\mu\text{-bd}}\to \theta$, that contradicts the assumption    $\lambda\rightarrow[\lambda]^{{<}\mu\text{-bd}}_{\theta}$. 

For $x\in [\lambda]^{{<\mu}\text{-bd}}$, let $\{\xi_i \mid i < \otp(x)\}$ be the increasing enumeration of $x$. define $F(x)$ to be $g(\delta_0, \dots, \delta_{n-1})$ where $\delta_i = f(\{\xi_j \mid j \in a_i\})$ and $h(\otp x) = \langle a_i \mid i < n\rangle$.

Let us claim that for every $A \in [\lambda]^\lambda$, $F\image [A]^{{<}\mu\text{-bd}} = \theta$. 

Let $B = f \image [A]^{{<}\mu\text{-bd}}$. By the assumption, $|B| = \theta$. By the choice of $g$, $g\image [B]^{{<}\omega} = \theta$. Let $\gamma \in \theta$. Let us pick $\langle \delta_0, \dots, \delta_{n-1}\rangle \in B^{<\omega}$ such that $g(\delta_0, \dots, \delta_{n-1}) = \gamma$. Let $t_i \in A^{{<}\mu\text{-bd}}$ be a set of ordinals such that $f(t_i) = \delta_i$. Let $y = t_0 \cup \dots \cup t_{n-1}$. Note that $y$ is still a bounded set of size $<\mu$. Let $\langle \zeta_i \mid i < \otp(y)\rangle$ be the increasing enumeration of the elements of $y$ and let $a_i = \{j < \otp(y) \mid \zeta_j \in t_j\}$. 

By the assumption on $h$ there is a ordinal $\rho \in [\otp(y), \mu)$ such that $h(\rho) = \langle a_0, \dots, a_{n-1}\rangle$. Let $x$ be $y \cup y'$ where $\otp(x) = \rho$, $\min y' > \sup y$, and $y' \subseteq A$. Clearly, $F(x) = \gamma$.  
\hfill \qedref{lemma:notJonsson}
\begin{theorem}
\label{thmmt2} Let $\lambda$ be I1. \newline 
Then one can force $\alpha_M(\lambda)=\mu^+$ where $\mu$ is a singular cardinal with uncountable cofinality.
\end{theorem}

\par\noindent\emph{Proof}. \newline 
Our starting point is a strongly Magidor cardinal $\lambda$ with a supercompact cardinal $\mu$ such that $\mu<\lambda$ and $\alpha_M^{<\mu}(\lambda)=\mu^+$. This can be arranged by part $(\beth)$ of Claim \ref{clmprepar}.

By Lemma \ref{lemma:notJonsson}, since $\mu^+$ is not J\'onsson, $\lambda \rightarrow [\lambda]^{<\mu\text{-bd}}_{\mu^+,<\mu^+}$. 

We shall force with Magidor forcing $\mathbb{M}$ to make $\mu>\cf(\mu)>\omega$ and our task is to show that $\alpha_M(\lambda)=\mu^+$ in the generic extension by $\mathbb{M}$.

Before proving this statement we need a preliminary assertion which reduces the number of the pertinent names for the coloring that we wish to define.
Assume that $\kappa$ is supercompact and $\mathbb{M}$ is Magidor forcing at $\kappa$. Let $A$ be a set of ordinals such that $|A|<\kappa$.
We claim that there exists a set $T$ of names, $|T|\leq\kappa$ such that for any name $\name{\tau}$ and any condition $p\in\mathbb{M}$ for which $p\Vdash \name{\tau}\subseteq\check{A} \wedge |\name{\tau}|=\aleph_0$ there exists a name $\name{\sigma}\in T$ and a condition $q\geq p$ so that $q\Vdash \name{\sigma} = \name{\tau}$.

For proving this assertion suppose that $|A|=\theta<\kappa$.
Fix a bijection $h:A\rightarrow\theta$.
Let $\mathcal{A} = \{p_i:i<\kappa\}$ be a maximal antichain of conditions which force an element into the Magidor sequence above $\theta$. 
For every $i<\kappa$ fix an ordinal $\alpha_i>\theta$ which is forced to be in the Magidor sequence by $p_i$.

Let $\langle d^i_1,\ldots,d^i_{n(i)}\rangle$ be the stem of $p_i$ for every $i<\kappa$.
There exists $m = m(i)\in[1,n(i)]$ such that $\alpha_i$ appears in $d^i_m$ (either $d^i_m=\alpha_i$ or $d^i_m=(\alpha_i,A_i)$).
A fundamental property of Magidor forcing is that $\mathbb{M}/p_i \cong \mathbb{M}_{\alpha_i}/p_i^{\leq m}\times \mathbb{M}/p_i^{>m}$ and new subsets of $\alpha_i$ are forced only by the lower part $\mathbb{M}_{\alpha_i}/p_i^{\leq m}$. The notation $\mathbb{M}_{\alpha_i}/p_i^{\leq m}$ should be understood as all conditions in the Magidor forcing below the cardinal $\alpha_i$, which are stronger than the condition $p_i^{\leq m}$. Notice that the number of names in $\mathbb{M}_{\alpha_i}/p_i^{\leq m}$ for subsets of $\alpha_i$ is at most $2^{\alpha_i}<\kappa$.
Let $T$ be the set of all names of the form $h^{-1}(y)$ where $y$ is a $\mathbb{M}_{\alpha_i}/p_i^{\leq m}$-name for a subset of $\theta$, for every $i<\kappa$, so $|T|\leq\kappa\cdot\kappa=\kappa$. 
We claim that $T$ is as required.

Indeed, assume that $\name{\tau}$ is an $\mathbb{M}$-name, $p\in\mathbb{M}$ and $p\Vdash \name{\tau}\subseteq\check{A} \wedge |\name{\tau}|=\aleph_0$.
By the maximality of $\mathcal{A}$ choose $p_i\in\mathcal{A}$ so that $p \parallel p_i$, and let $q\in\mathbb{M}$ be a condition which satisfies $q\geq p,p_i$. Now $q \Vdash \name{\tau}\subseteq\check{A} \wedge |\name{\tau}|=\aleph_0$ as $q\geq p$ and $q\Vdash \name{\sigma}=\name{\tau}$ for some $\name{\sigma}\in T$ since $q\geq p_i$.
This completes the proof of the assertion.

Let $G\subseteq\mathbb{M}$ be generic.
We try to show that $V[G]\models \alpha_M(\lambda)=\mu^+$.
By the fact that $\lambda$ is I1 in $V$ we can see that $\lambda$ is I1 (and hence Magidor) in $V[G]$.
Assume that $\name{f}:[\lambda]^{\aleph_0\text{-bd}} \rightarrow \mu^+$.
Let $A$ be a bounded subset of $\lambda$ of size less than $\mu$ which belongs to $V$. Define $g(A)$ to be the supremum of all ordinals $\alpha < \mu^+$ such that there is an $\mathbb{M}$-name $\name{\tau}$, the weakest condition forces that $\name{\tau} \subseteq \check{A} \wedge |\name{\tau}|=\aleph_0$, and there is a condition $q$ such that $q\Vdash \name{f}(\name{\tau}) = \check{\alpha}$.

By the preliminary assertion and the chain condition of $\mathbb{M}$ we see that $g(A)\in\mu^+$ for every $A$ as above.
Since $g\in V$ and $\alpha_M^{<\mu}(\lambda)=\mu^+$ there are $H\in[\lambda]^\lambda$ and $\beta\in\mu^+$ such that $g\image H^{<\mu\text{-bd}}\subseteq\beta$.
We can finish the proof by showing that $\name{f}\image[H]^{\aleph_0\text{-bd}}$ is forced to be a subset of $\beta$ as well.

Suppose not.
Choose $\gamma, \name{\tau}$ and $q$ such that $\gamma>\beta, q\in\mathbb{M}, \name{\tau}$ is an $\mathbb{M}$-name of a countable bounded subset of $H$ and $q\Vdash\name{f}(\name{\tau})=\gamma$.
Fix $A\in V,\, A\subseteq H,\,|A|=\theta<\mu$ and $r\geq q$ such that $r\Vdash\name{\tau} \subseteq A$.
Notice that $r\Vdash\name{\tau}\subseteq H$ as well (since $r\geq q$) and hence for some condition $s\geq r$ we have $s\Vdash \name{f}(\name{\tau})\leq g(A)<\beta<\gamma$, a contradiction.

We showed that in the generic extension, $\alpha_M(\lambda) \leq \mu^+$. Note that $\alpha_M(\lambda) \neq \mu$ since $\mu$ is singular. The set of all $\zeta < \mu$ in which $\mathbb{M}$ adds a Prikry sequence is unbounded in $\mu$. Let $\zeta$ be a measurable cardinal of Mitchell order $1$ in the generic Magidor club. Let $p\in\mathbb{M}$ be a condition that forces $\zeta$ to be in the Magidor club. The forcing $\mathbb{M} / p$ decomposes into a product of two forcing notions, $\mathbb{M}' \times \mathbb{P}$ where $\mathbb{P}$ is the Prikry forcing at $\zeta$. By standard arguments, one can verify that $\mathbb{P}$ is equivalent to the Prikry forcing at $\zeta$ in the generic extension by $\mathbb{M}'$. Thus, $\mathbb{P}$ forces that $\alpha_M(\lambda) > \zeta$. 

Since this is true for all $\zeta < \mu$, $\alpha_M(\lambda) \geq \mu^+$ and thus $\alpha_M(\lambda) = \mu^+$, as wanted. 

\hfill \qedref{thmmt2}

As in the former section, the above proof can be abstracted. Basically all we need is the special covering property of new countable sets by relatively small sets. This is the main distinction between Prikry and Magidor forcing, used above.

We conclude the paper with an open problem which goes back to singular cardinals with countable cofinality.
The main theorem of this section states that $\alpha_M$ can realize the true cofinality of a product of cardinals below some $\mu>\cf(\mu)>\omega$. 
The following is natural:

\begin{question}
\label{qtcf} Assume that $\mu<\lambda$, $\lambda$ is Magidor and $\mu>\cf(\mu)=\omega$. \newline 
Is it possible that $\alpha_M = {\rm tcf}(\prod_{n\in\omega}\mu_n,J^{\rm bd}_\omega)$ for some increasing sequence of regular cardinals $\langle \mu_n:n\in\omega\rangle$?
\end{question}

\section{Acknowledgments}
First author's research is supported by Shelah's ERC grant 338821.

We would like to thank the anonymous referee for many helpful suggestions that improved this paper's readability and correctness.
\newpage

\providecommand{\bysame}{\leavevmode\hbox to3em{\hrulefill}\thinspace}
\providecommand{\MR}{\relax\ifhmode\unskip\space\fi MR }
\providecommand{\MRhref}[2]{%
  \href{http://www.ams.org/mathscinet-getitem?mr=#1}{#2}
}
\providecommand{\href}[2]{#2}


\end{document}